\documentclass[12pt]{article}
\usepackage{arxiv}

\usepackage{siunitx}
\usepackage[utf8]{inputenc} 
\usepackage[T1]{fontenc}    
\usepackage{amsmath}
\usepackage{bm}
\usepackage{booktabs}       
\usepackage{amsfonts}       
\usepackage{amssymb}
\usepackage{nicefrac}       
\usepackage{microtype}      
\usepackage{soul}
\usepackage[table]{xcolor}
\definecolor{medium_blue}{rgb}{0, 0, 0.804}
\usepackage{hyperref}
\hypersetup{colorlinks=true,linkcolor=medium_blue, citecolor=medium_blue}
\usepackage{leftidx}
\usepackage{graphicx}
\usepackage{threeparttable}
\usepackage{multirow}
\usepackage{float}
\usepackage{physics}
\usepackage{lscape}
\usepackage{threeparttable}

\newcommand{\RomanNumber}[1]{\MakeUppercase{\romannumeral #1}}

\title{A dual adaptive explicit time integration algorithm for efficiently solving the cardiac monodomain equation}

\author{
  Konstantinos A. ~Mountris\thanks{[mail] kmountris@unizar.es \quad [url] https://www.mountris.org} \\
  Arag\'on Institute of Engineering Research, IIS Arag\'on\\
  CIBER-BBN\\
  University of Zaragoza\\
  Spain, Zaragoza, ZGZ 50018 \\
  \texttt{kmountirs@unizar.es} \\
  \And
  Esther ~Pueyo \\
  Arag\'on Institute of Engineering Research, IIS Arag\'on, \\
  CIBER-BBN\\
  University of Zaragoza\\
  Spain, Zaragoza, ZGZ 50018 \\
  \texttt{epueyo@unizar.es} \\
}

\begin{document}
\maketitle

\begin{abstract}
The monodomain model is widely used in \textit{in silico} cardiology to describe excitation propagation in the myocardium. Frequently, operator splitting is used to decouple the stiff reaction term and the diffusion term in the monodomain model so that they can be solved separately. Commonly, the diffusion term is solved implicitly with a large time step while the reaction term is solved by using an explicit method with adaptive time stepping. In this work, we propose a fully explicit method for the solution of the decoupled monodomain model. In contrast to semi-implicit methods, fully explicit methods present lower memory footprint and higher scalability. However, such methods are only conditionally stable. We overcome the conditional stability limitation by proposing a dual adaptive explicit method in which adaptive time integration is applied for the solution of both the reaction and diffusion terms. In a set of numerical examples where cardiac propagation is simulated under physiological and pathophysiological conditions, results show that our proposed method presents preserved accuracy and improved computational efficiency as compared to standard operator splitting-based methods.
\end{abstract}


\section{Introduction} \label{sec:introduction}
Computational modeling and simulation is widely used in cardiac electrophysiology to gain more insight into the mechanisms underlying the heart's electrical activity, predict unfavorable responses in the presence of disease or identify novel therapeutic targets. Propagation of action potential (AP) waves in the myocardium can be simulated by solving a system of partial differential equations (PDEs) known as the bidomain model \cite{tung1978}. The bidomain model considers the cardiac tissue as a continuum of two anisotropic compartments describing the intracellular and extracellular spaces. Assuming equal anisotropy ratios for the intracellular and extracellular spaces, the bidomain model can be reduced to the simplified monodomain model \cite{keener2009}. In situations where the currents in the extracellular domain have little influence on cardiac transmembrane potential and ionic currents, the monodomain model can produce realistic activation patterns and transmembrane potential values with less computational cost than the bidomain model \cite{potse2006}.

\noindent The monodomain model is described by a single reaction-diffusion PDE for the transmembrane potential across the myocardium, while the extracellular potential can be computed from another PDE once the transmembrane potential has been solved. In the reaction-diffusion PDE, the reaction term describes the generation of the cellular AP and the diffusion term describes its propagation in the tissue. Realistic cardiac cell models are composed of a large set of stiff ordinary differential equations (ODEs) to describe the temporal evolution of the ionic concentrations and gating variables of the cell. Solving stiff ODE models can be time consuming, especially for large scale problems, since a small time integration step is required to ensure numerical stability. The operator splitting technique \cite{qu1999} can be used to decouple the stiff reaction term and the diffusion term. In this way, a larger time step can be used to solve the diffusion term independently of the reaction term, which may be solved adaptively using a small time step \cite{gomes2020,heidenreich2010}.

\noindent Commonly, the decoupled reaction-diffusion system is solved by employing a semi-implicit scheme. The stiff reaction term is integrated using an explicit time integration method (e.g. forward Euler, Runge-Kutta 4) whereas the diffusion term is solved using an implicit time integration method (e.g. backward Euler, Crank-Nicolson). Implicit time integration methods are popular for their unconditional stability. However, they require solving a system of equations at each time step, which makes them more complex to implement and harder to parallelize than explicit methods. Explicit methods, on the other hand, are only conditionally stable. Yet attractive, explicit methods are impractical for large scale problems where a high mesh resolution is required. This is so because a very small time step is required to ensure stability, as with decreasing mesh spacing the upper bound of the time step is decreased too \cite{courant1967}.

\noindent In this work, we identify a simple, yet meaningful, realization that allows overcoming the conditional stability limitation of the explicit time integration scheme. Once the reaction and diffusion terms of the monodomain model are decoupled by application of the operator splitting technique, we propose a dual adaptive explicit time integration method where both the reaction and diffusion terms are solved explicitly with a different adaptive scheme in each case. Since most of the computational burden is associated with the solution of the reaction term, the overhead of the adaptive solution of the diffusion term is minimum. The structure of the article is as follows. In section \ref{sec:dual_adaptive}, we describe the dual adaptive explicit time integration algorithm. In section \ref{sec:eval_dae}, we compare the accuracy and efficiency of the proposed method against the standard method with adaptive time integration for the reaction term only and against the method without any time step adaptation. The comparison is performed for two-dimensional (2D) and three-dimensional (3D) problems of cardiac electrophysiology in both health and disease. In section \ref{sec:discuss}, we discuss our findings and in section \ref{sec:conclude}, we present the concluding remarks.

\section{Dual adaptive explicit time integration} \label{sec:dual_adaptive}

Propagation of the cardiac AP was simulated by using the monodomain model given by:
\begin{equation} \label{eq:monodomain}
\begin{array}{ll}
    \partial V / \partial t = -I_{ion}(V) / C + \bm{\nabla} \cdot (\bm{D} \bm{\nabla} V) &\textrm{ in } \Omega \\
    \bm{n} \cdot (\bm{D}\bm{\nabla}V) = 0 &\textrm{ in } \partial \Omega
\end{array}
\end{equation}
\noindent where $\partial V / \partial t$ is the time derivative of the transmembrane potential, $I_{ion}$ is the total ionic current, $C$ denotes the cell capacitance per unit surface area and $\bm{D}$ is the diffusion tensor. $\Omega$ and $\partial \Omega$ denote the domain of interest and its boundary, respectively, and $\bm{n}$ is the outward unit vector normal to $\partial \Omega$. 

\noindent Employing the operator splitting method \cite{qu1999}, equation (\ref{eq:monodomain}) can be written as:
\begin{subequations}  \label{eq:opsplit}
\begin{gather}
    \partial V / \partial t = -I_{ion}(V) / C \textrm{ in } \Omega \label{eq:opsplit_a} \\
    \partial V / \partial t = \bm{\nabla} \cdot (\bm{D} \bm{\nabla} V) \textrm{ in } \Omega \label{eq:opsplit_b} \\
    \bm{n} \cdot (\bm{D}\bm{\nabla}V) = 0 \textrm{ in } \partial \Omega \label{eq:opsplit_c}
\end{gather}
\end{subequations}
Following Strang's operator splitting, time integration of the set of equations (\ref{eq:opsplit}) for the time interval $\left [0, T\right]$ using a time step $dt$ is performed in three steps (\RomanNumber{1} - \RomanNumber{3}). Taking the solution at a given time $t$ as the initial condition, adaptive integration is applied in step \RomanNumber{1} to solve the diffusion term (equations \ref{eq:opsplit_b} and \ref{eq:opsplit_c}) with time step $dt/2$. Using the results of step \RomanNumber{1} as initial condition, adaptive integration for the reaction term (equation \ref{eq:opsplit_a} and ODEs for ionic concentrations and gating variables used in the computation of $I_{ion}(V)$) is applied at step \RomanNumber{2} for time step $dt$. Using the results of \RomanNumber{2}, the iteration terminates by integrating again the diffusion term for time step $dt/2$ in step \RomanNumber{3}. In practice, steps \RomanNumber{1} and \RomanNumber{3} can be combined into only one step (denoted as step B), except for the initial and final steps of the integration in the interval $\left [0, T\right]$. Step \RomanNumber{2} is denoted as step A. 

\noindent In our proposed dual adaptive explicit time integration (DAETI) method, similarly to the standard adaptive time integration method \cite{qu1999}, the reaction term integration in step \RomanNumber{2} was performed by using the forward Euler method with adaptive time step $dt_{\textrm{ar}} = dt/k$, where the parameter $k$ was defined to be an integer so as to keep steps A and B, corresponding to the reaction and diffusion terms, synchronized every time step $dt$. The value of $k$ was selected in the range $\left[1, k_{\textrm{max}}\right]$. The upperbound $k_{\textrm{max}}$ was obtained as $ k_{max} = \left \lfloor dt/dt_0 \right \rfloor $, where $dt_0$ is the maximum value of the time step that guarantees numerical stability of the cell electrophysiology model. In each iteration, the value of the parameter $k$ was calculated by $k = k_0 + \lfloor \abs{\partial V/\partial t} \rfloor$, where $\partial V/\partial t$ is the time derivative of the AP at the previous iteration. For $\partial V/\partial t > 0$ (e.g. steep gradient during upstroke), we chose $k_0 = 5$ to ensure safe propagation of the wave front. For $\partial V/\partial t \leq 0$ (e.g. smooth gradient during repolarization), we chose $k_0 = 1$ to avoid zero-division during the calculation of $dt_{\textrm{ar}}$. If $k > k_{\textrm{max}}$, then $k = k_{\textrm{max}}$. For a more detailed explanation on the choice of $k_0$, we refer the reader to \cite{qu1999}.

\noindent The adaptive time integration of the diffusion term in steps \RomanNumber{1} and \RomanNumber{3} was performed using the forward Euler method with time step $dt_{\textrm{ad}} = dt/2l$, where $l = \lfloor dt / 2dt_{\textrm{s}} \rfloor$ if $dt/2 > dt_{\textrm{s}}$ and $l = 1$, otherwise. A stable diffusion time step $dt_{\textrm{s}}$ was obtained by the Gerschg\"{o}rin theorem \cite{myers1978}:  
\begin{equation} \label{eq:gerschgorin}
    dt_{\textrm{s}} = 0.9  \min_{i=1,...,n} \left[ \frac{m_{ii}}{k_{ii} + \sum\limits_{\substack{j=1 \\ j \neq i}}^{n} \abs{k_{ij}}} \right].
\end{equation}
\noindent where $m_{ij}$, $k_{ij}$ are the elements of the $i$-th row and $j$-th column of the mass and stiffness matrices, respectively, for the Finite Element approximation of equation (\ref{eq:opsplit_b}). The multiplication factor 0.9 is a safety factor to ensure numerical stability for the estimated time step. When $dt/2 > dt_{\textrm{s}}$, the adaptive integration of the diffusion term ensures numerical stability, while when $dt/2 \leq dt_{\textrm{s}}$, the DAETI method is reduced to the standard operator splitting-based integration.

\section{Evaluation of dual adaptive explicit time integration} \label{sec:eval_dae}

\noindent To evaluate the accuracy and efficiency of the DAETI method we performed simulations under physiological and pathophysiological conditions for a 2D tissue and a 3D biventricular geometry. We compared the simulation results obtained by DAETI with simulation results obtained by the operator splitting technique with no adaptive integration (OST) and by OST with adaptive integration of the reaction term (OSTAR). Both the reaction and diffusion terms in OST and OSTAR were solved explicitly using the forward Euler method as in DAETI. In all simulations, the time integration step was $dt = 0.1$ ms for DAETI and $dt = 0.01$ ms for OST. In OSTAR, $dt = 0.1$ ms was used for coarse mesh discretizations where $dt_{\textrm{s}} \geq dt$. For fine mesh discretizations where $dt_{\textrm{s}} < dt$, $dt_{\textrm{s}}$ was used as the time integration step. Simulations were performed using a multithreaded implementation of the Finite Element Method using ELECTRA, an in-house software implementing the Finite Element Method and the Meshfree Mixed Collocation method \cite{mountris2019,mountris2020cinc,mountris2020radial} for solving the monodomain model. In this work, we used the Finite Element implementation. All simulations were performed on a laptop with Intel\textsuperscript{\textregistered} Core\texttrademark i7-4720HQ CPU and 16 GB of RAM.

\subsection{Electrical propagation in a 2D cardiac tissue} \label{subsec:ap_2D}

AP propagation in a $5\cross5$ cm ventricular tissue sheet was simulated for $T=500$ ms after achieving steady-state conditions at a pacing cycle length of 1000 ms. A stimulus of 1 ms duration and twice diastolic threshold amplitude was delivered at the left side of the tissue ($X=0$ cm) at time $t = 50$ ms. 4-node regular isoparametric elements with a spacing step of $h = 100$ $\mu$m were used. A homogeneous ventricular tissue was considered, with cellular electrophysiology represented by the O{'}Hara et al. human epicardial AP model \cite{ohara2011}. Fiber orientation was considered parallel to the X-axis. The tissue diffusion coefficient was $d_0 = 0.0017$ cm$^2$/ms with transverse-to-longitudinal ratio $\rho = 0.25$. The critical value for the diffusion time step $dt_{\textrm{s}}$ was estimated using equation (\ref{eq:gerschgorin}), rendering a value of $dt_{\textrm{s}} = 0.013$ ms.

\begin{figure}[htbp]
\centering
\includegraphics[width=8cm]{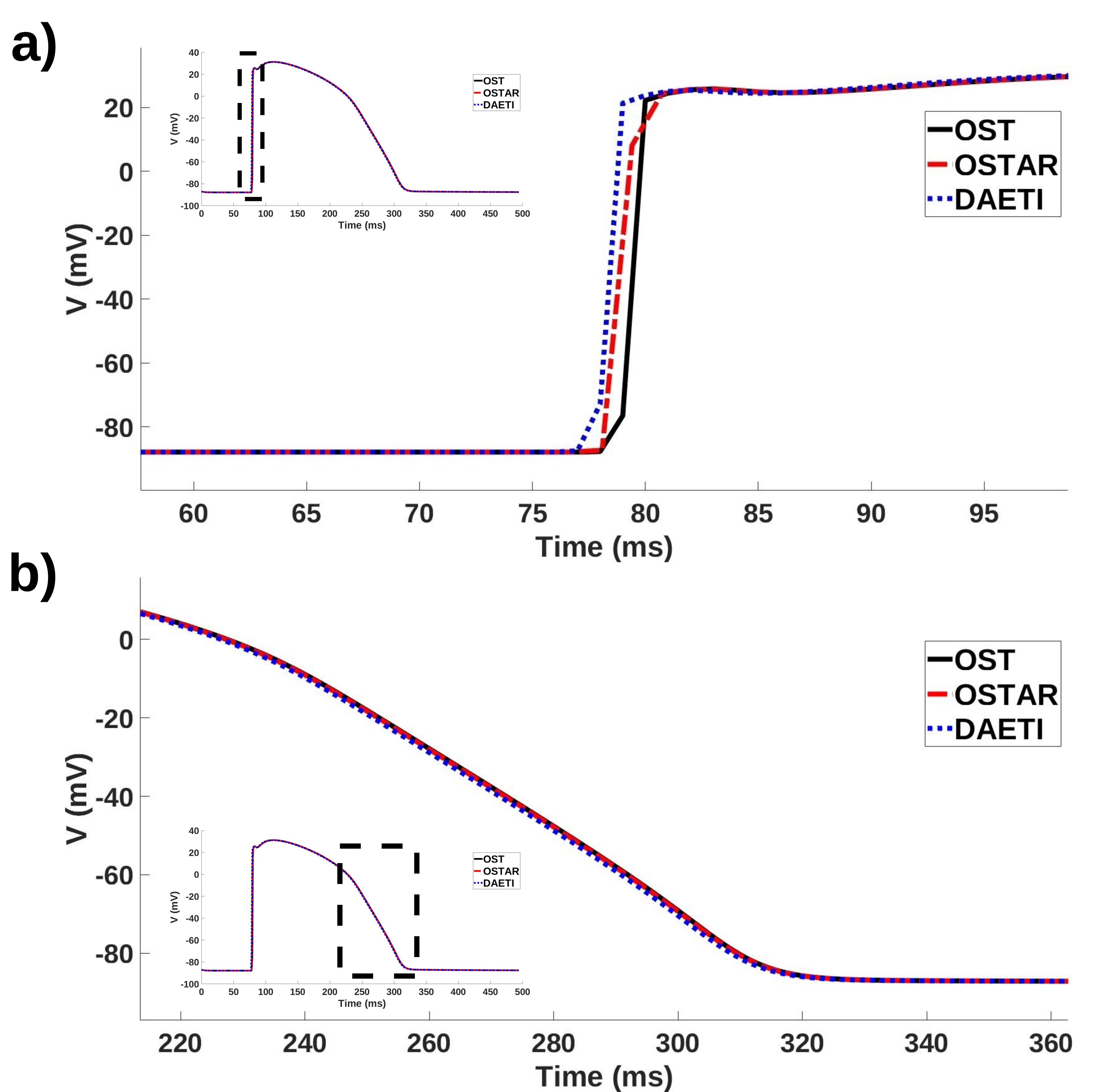}
\caption{APs at the center of a $5\cross5$ cm human ventricular epicardial tissue calculated using OST, OSTAR and DAETI methods. a) Comparison of APs at upstroke. b) Comparison of APs at repolarization.}
\label{fig:ap_comparison}
\end{figure}

\noindent The APs at the center of the tissue sheet calculated using the three evaluated methods are shown in Figure \ref{fig:ap_comparison}. As can be observed from the figure, differences in transmembrane potential between DAETI and any of the other two methods were minimal, being within 0.34 mV in all cases for aligned APs. The longitudinal conduction velocity (CV) obtained by using the DAETI method was CV = 0.0633 cm/ms, while CV values for OSTAR and OST methods were CV = 0.0614 cm/ms and CV = 0.0595 cm/ms, respectively. AP duration at 90\% repolarization (APD$_{90}$) was APD$_{90}$ $= 258$ ms for the three simulation methods. The total execution time was 21 min for DAETI, 95 min for OSTAR, and 103 min for OST. 

\subsection{Electrical propagation in a 2D cardiac tissue with fibrosis} \label{subsec:cross-field}

We next simulated AP propagation in a $5 \cross 5$ cm ventricular epicardial tissue sheet, as in subsection \ref{subsec:ap_2D}, but in this case we included $10\%$ diffusive fibrosis by randomly distributing fibroblasts across the tissue following a uniform distribution. Human ventricular epicardial cell electrophysiology was represented by the O\'Hara et al. model, while fibroblast electrophysiology was represented by the MacCannell active fibroblast model \cite{maccannell2007}. The tissue was stimulated by applying a cross-field stimulation protocol to generate a sustained spiral wave. After achieving steady-state at a pacing cycle length of 1000 ms, a S1 stimulus was delivered at the left side of the tissue ($X=0$ cm) at time $t = 50$ ms, which was followed by a premature S2 stimulus delivered at the bottom left corner of the tissue ($X=0-1.25$ cm, $Y=0-2.50$ cm) at time $t=200$ ms. 

We performed simulations on six meshes with 4-node regular elements and varying spacing steps ranging from $h = 100$ to $h = 200$ $\mu$m. The computational efficiency of DAETI was compared against that of OSTAR and OST. The critical diffusion time step $dt_{\textrm{s}}$ was obtained by equation (\ref{eq:gerschgorin}) for each mesh. The diffusion coefficient was $d_{0m} = 0.002$ cm$^2$/ms between epicardial myocytes and $d_{0f} = 0.00066$ cm$^2$/ms between fibroblasts as well as for the interaction between myocytes and fibroblasts. The transverse-to-longitudinal conductivity ratio was set to $\rho = 0.25$. The characteristics of the simulated meshes and the corresponding values of $dt_{\textrm{s}}$ are provided in Table \ref{tab:statistics}.

\begin{table}[htbp]
\centering
\caption{Summary of mesh characteristics used in the simulation of a 2D fibrotic tissue under cross-field stimulation. $h$ denotes the mesh spacing and $dt_{\textrm{s}}$, the critical value for the diffusion time integration step.}
\label{tab:statistics}
\begin{tabular}{cccc}
\hline
$h$ ($\mu$m) & Nodes No. & Elements No. & $dt_{\textrm{s}}$ (ms) \\
\hline
\hline
200 & 63001 & 62500 & 0.113 \\
180 & 77284 & 76729 & 0.091 \\
160 & 97969 & 97344 & 0.072 \\
140 & 128164 & 127449 & 0.059 \\
120 & 173889 & 173056 & 0.044 \\
100 & 251001 & 250000 & 0.028 \\
\hline
\end{tabular}
\end{table}

\begin{figure}[htbp]
\centering
\includegraphics[width=9cm]{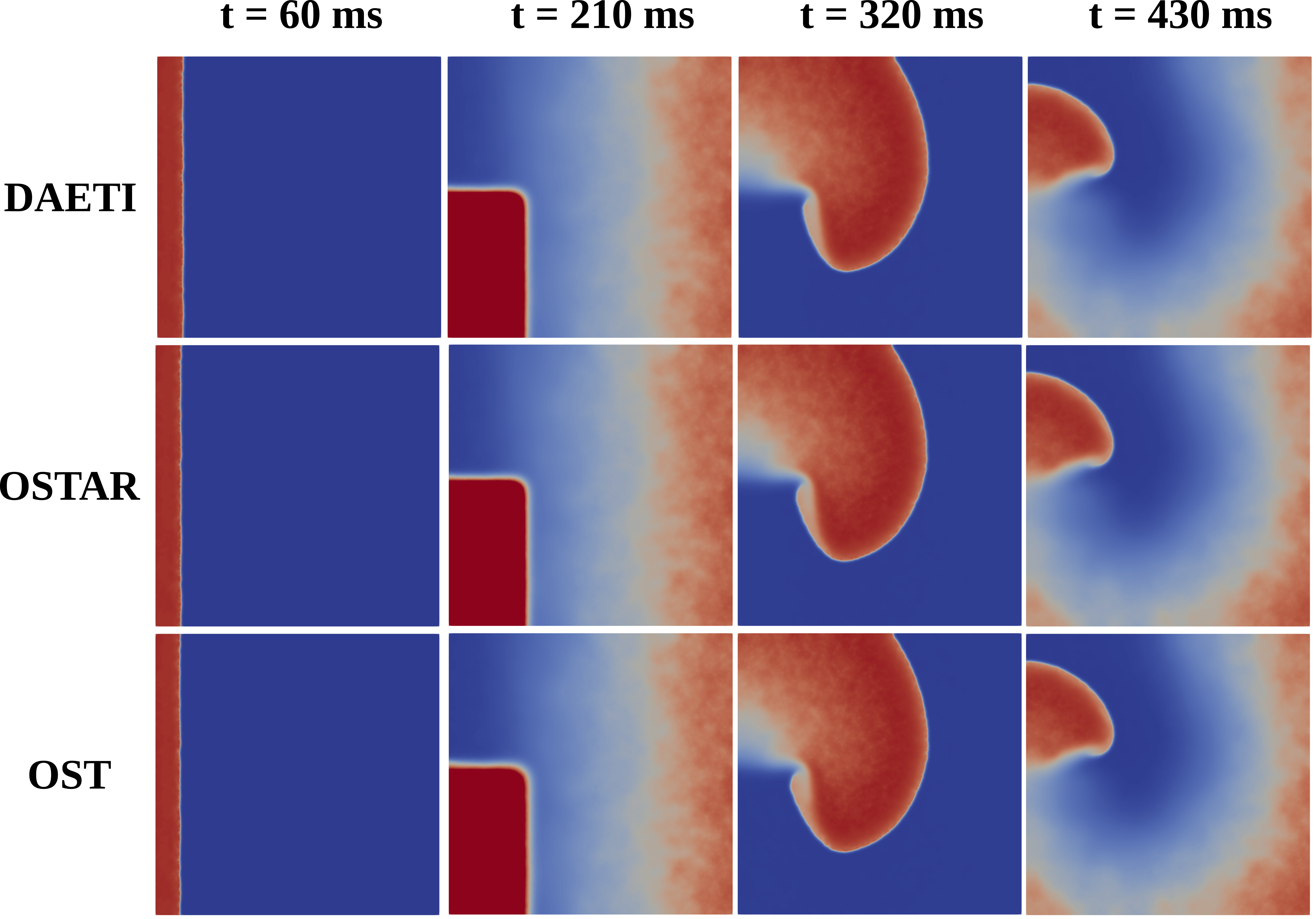}
\caption{Voltage snapshots in a 2D fibrotic tissue sheet in response to a cross-stimulation protocol. Spiral waves of the same characteristics can be observed for DAETI (top), OSTAR (middle) and OST (bottom) at different simulated time instants.}
\label{fig:rotors}
\end{figure}

Transmembrane voltage in the simulated tissues after application of the cross-field stimulation protocol is presented in Figure \ref{fig:rotors} for the mesh with $h=100$ $\mu$m. The three methods rendered highly similar voltage values, being the differences between aligned APs below 0.41 mV, with the three of them presenting the same characteristics of the generated spiral waves at the different time instants along the simulation time. 

\begin{figure}[htbp]
\centering
\includegraphics[width=9cm]{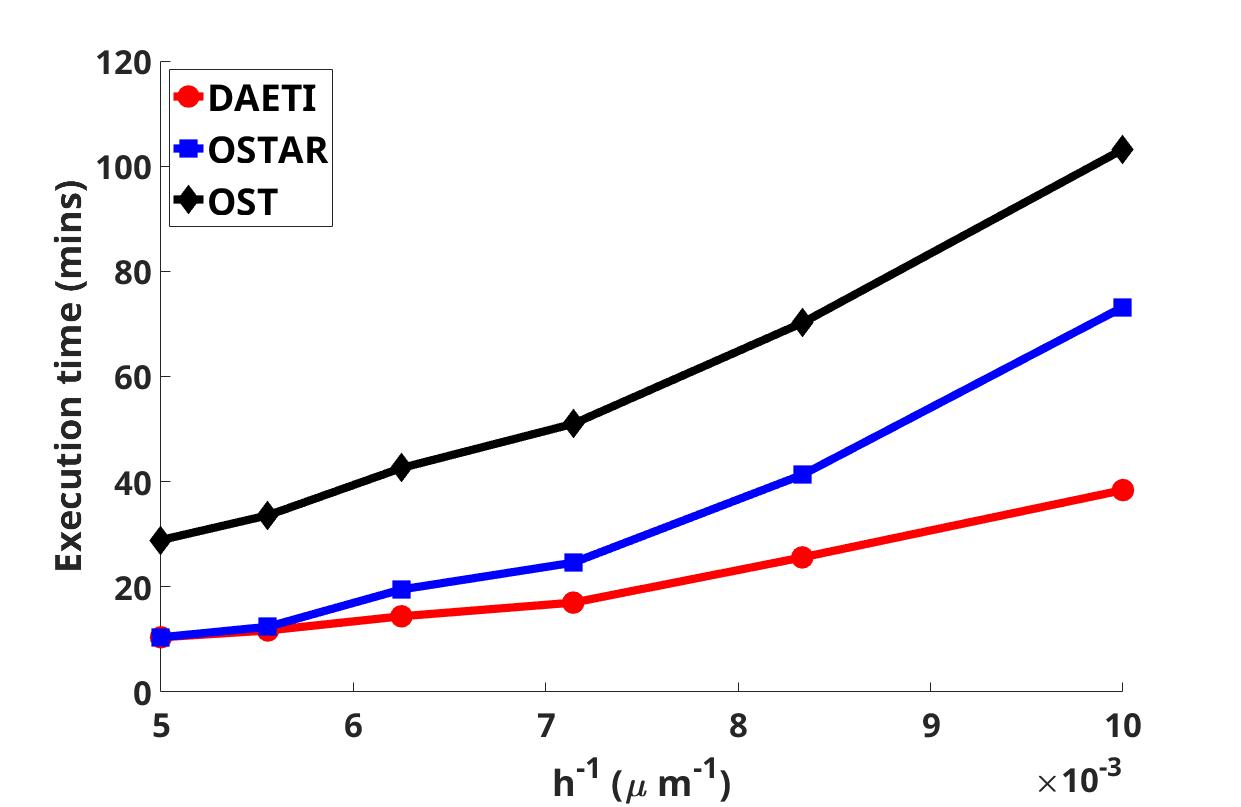}
\caption{Execution time for DAETI (red circle), OSTAR (blue square) and OST (black diamond) evaluated in a 2D fibrotic tissue with 4-node regular element meshes and varying space steps ranging from $h = 100$ to $h = 200$ $\mu$m.}
\label{fig:time_comparison}
\end{figure}

The total execution time for a simulation time of $T=500$ ms is given in Figure \ref{fig:time_comparison} for the three evaluated methods. As can be observed from the figure, the computational time required by DAETI was notably lower than that of OST. When compared to OSTAR, DAETI was associated with similar computational times for coarse meshes but with remarkably lower times for fine meshes.

\subsection{Electrical propagation in a 3D cardiac biventricular geometry} \label{sec:biventricular}
We further compared the DAETI method with the OSTAR and OST methods in the simulation of AP propagation in a 3D biventricular swine geometry at baseline and under left bundle branch block (LBBB) conditions. A tetrahedral mesh representation of the biventricular swine anatomy (273919 nodes, 1334218 elements) was made available from the CRT-EPiggy19 challenge \cite{camara2019,pop2020}. It was part of a swine model dataset of LBBB for experimental studies of cardiac resynchronization therapy \cite{rigol2013,iglesias2016}. A homogeneous diffusion coefficient $d_0 = 0.002$ cm$^2$/ms was set across the cardiac tissue, with a transverse-to-longitudinal conductivity ratio of $\rho = 0.25$. The diffusion coefficient for the connective tissue at the base of the biventricular anatomy was set to $d_0 = 0.00066$ cm$^2$/ms. The electrophysiology of connective tissue was represented by the MacCannell active fibroblast model. For the rest of the biventricular geometry, the O{'}Hara et al. cell model was used. Transmural heterogeneities were included by considering endocardial, midmyocardial and epicardial cells across the ventricular wall at 0.5:0.2:0.3 ratio. The orientation of the myocardial fibers was computed using a rule-based method \cite{doste2019}.

\noindent A Purkinje conduction system was generated for the biventricular model by using a fractal-tree generation algorithm \cite{costabal2016}. Purkinje-myocyte junctions (PMJs) on the endocardial surface were obtained by applying a range-search algorithm using a spherical search area with radius $R=2$ mm, centered at each end node of the conduction system. Rather than considering AP propagation in the conduction system, periodic stimuli of 1 ms duration and twice diastolic threshold in amplitude were applied to each PMJ at a cycle length of 1000 ms. To ensure realistic activation at baseline conditions \cite{durrer1970}, we divided PMJs into four groups: left apex - LA; left base - LB; right apex - RA; right base - RB. LA-PMJs were activated at time $t = 0$ ms, while the activation of LB-PMJs, RA-PMJs, and RB-PMJs was delayed by 7, 4, and 11 ms, respectively. Under LBBB conditions, stimulation at LA-PMJs and LB-PMJs was blocked. After achieving steady-state, the total simulation time was $T = 500$ ms for both baseline and LBBB. The critical value for the diffusion time step was $dt_{\textrm{s}} = 0.034$ ms.

\noindent Baseline and LBBB simulation results obtained by DAETI and OSTAR were compared against the results obtained by the OST simulation, taken here as a reference, by computing the normalized root mean square error for both the local activation time (LAT), denoted by $e_{\textrm{LAT}}$, and APD$_{90}$, denoted by $e_{\textrm{APD}}$:
\begin{equation} \label{eq:nrmse}
    e_{\textrm{Z}} = \frac{\displaystyle \sqrt{\frac{1}{N}\sum_{i=1}^N (\hat{u}_i - u_i)^2}} {u_{max} - u_{min}},
\end{equation}
\noindent where the subindex Z stands for LAT or APD, $N$ is the number of nodes in the mesh, $u$ denotes the reference nodal value for either LAT or APD$_{90}$ obtained by OST and $\hat{u}$ denotes the nodal LAT or APD$_{90}$ value obtained by OSTAR or DAETI. $u_{max}$ and $u_{min}$ denote the maximum and minimum values of $u$ across all nodes in the tissue.

\noindent Using the DAETI method, $e_{\textrm{LAT}}$ was 3.5E$^{-3}$ and 4.4E$^{-3}$ at baseline and under LBBB conditions, respectively, whereas for the OSTAR method the corresponding $e_{\textrm{LAT}}$ values were 2.4E$^{-3}$ and 2.6E$^{-3}$. In the case of $e_{\textrm{APD}}$, the values for DAETI were 8.3E$^{-4}$ and 7.8E$^{-4}$ at baseline and under LBBB conditions, while for OSTAR these were 6.7E$^{-4}$ and 7.1E$^{-4}$. 

\begin{figure}[htbp]
\centering
\includegraphics[width=10cm]{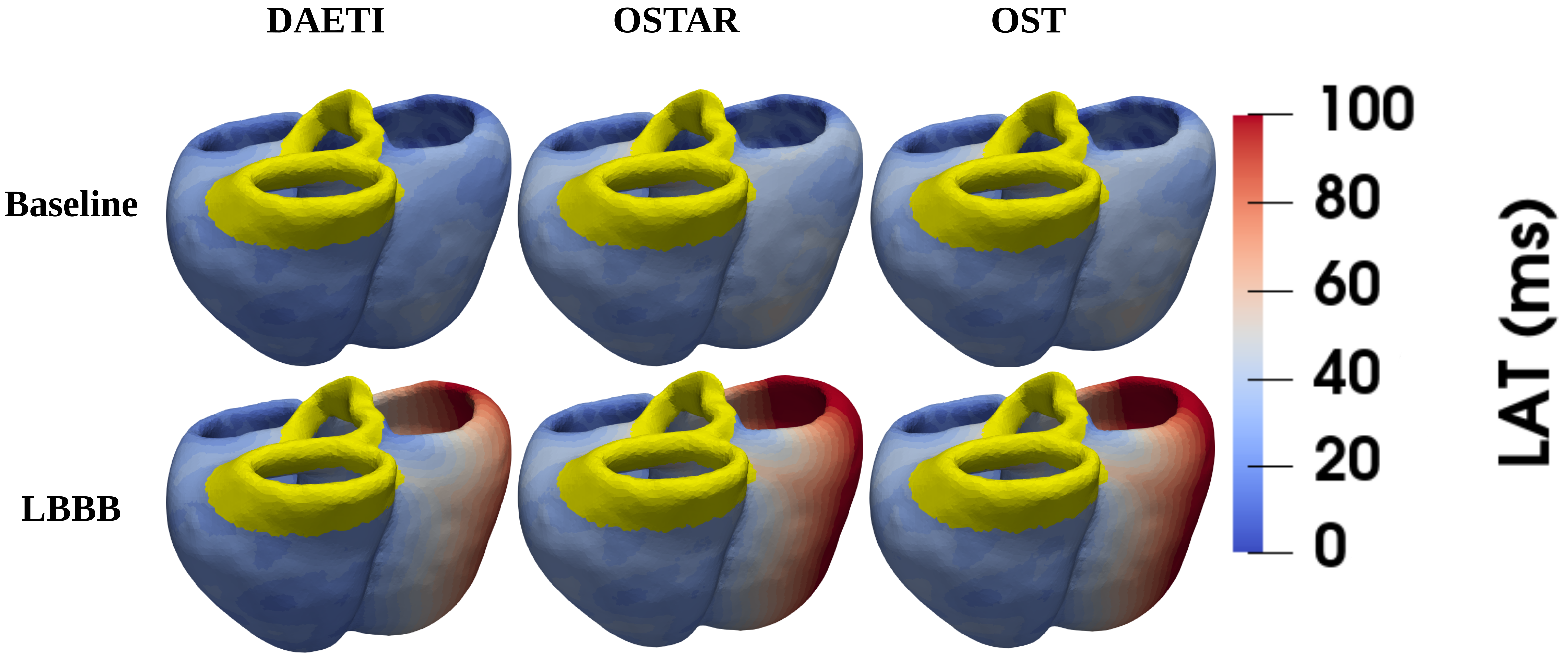}
\caption{LAT maps at baseline (top) and under LBBB conditions (bottom) for a simulation using DAETI (left), OSTAR (middle) and OST (right) methods.}
\label{fig:lat3D}
\end{figure}

\begin{figure}[htbp]
\centering
\includegraphics[width=10cm]{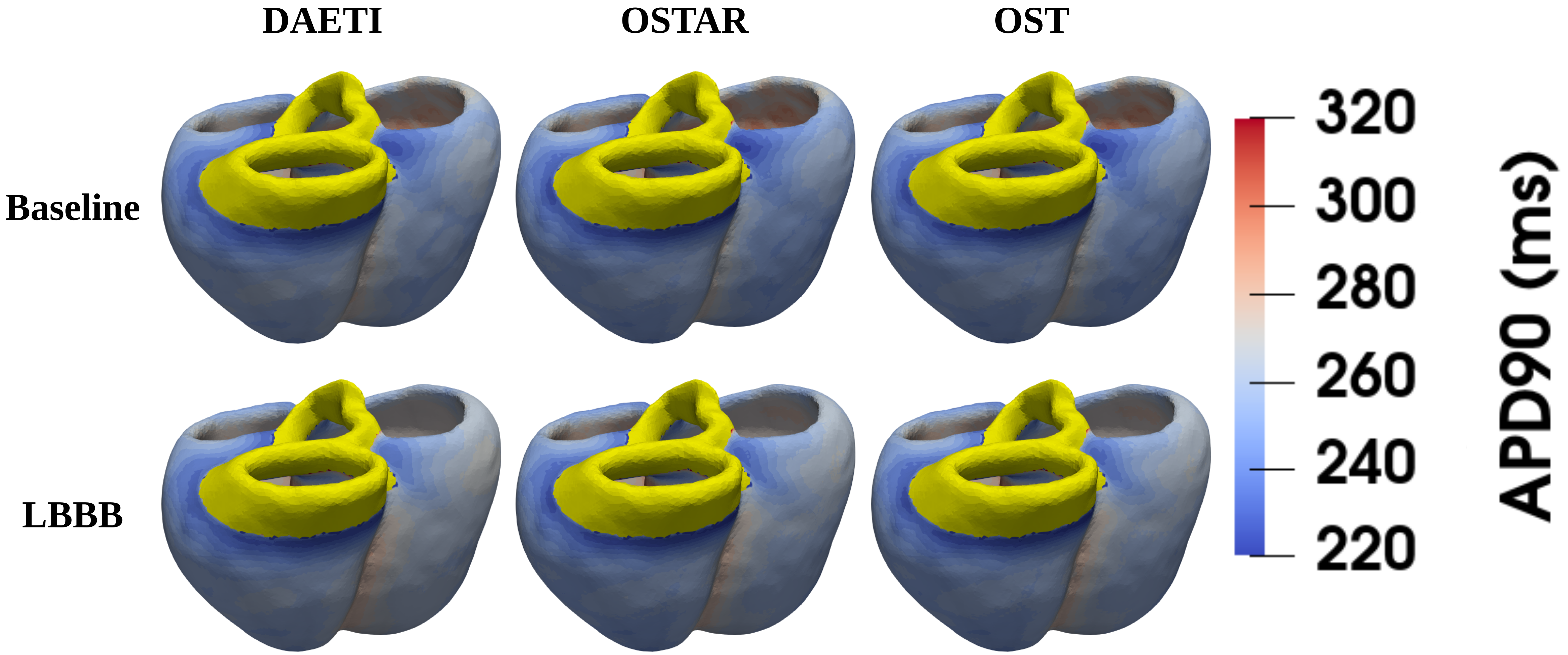}
\caption{APD$_{90}$ maps at baseline (top) and under LBBB conditions (bottom) for a simulation using DAETI (left), OSTAR (middle) and OST (right) methods.}
\label{fig:apd3D}
\end{figure}

\noindent Mean LAT was 23.8 ms, 24.1 ms, and 24.2 ms for DAETI, OSTAR, and OST at baseline. Under LBBB, mean LAT was 49.3 ms, 49.7 ms, and 49.9 ms, respectively. Mean APD$_{90}$ was 233.5 ms at baseline and 233.2 ms under LBBB for all three methods. Figures \ref{fig:lat3D} and \ref{fig:apd3D} show LAT and APD$_{90}$ maps at baseline and under LBBB. From Figure \ref{fig:lat3D}, it can be appreciated that the epicardial activation under LBBB occurs by wave propagation through the interventricular septum and the anterior and posterior left ventricular wall. This characteristic activation pattern, as well as the maximum LAT for LBBB being 120 ms larger than at baseline, are well reproduced by the three numerical integration methods.

\noindent In terms of computational efficiency, the total execution time at baseline was 25.6 min for DAETI, 60.6 min for OSTAR, and 106.8 min for OST. Under LBBB conditions, the total execution time was 23.4 min for DAETI, 62.4 min for OSTAR, and 102.6 min for OST. Using the DAETI method we obtained a speed-up of x4.1 at baseline and x4.4 under LBBB with respect to OST. The speed-up obtained by OSTAR with respect to OST was x1.8 at baseline and x1.7 under LBBB.

\section{Discussion} \label{sec:discuss}

We proposed a dual adaptive explicit time integration (DAETI) method to solve the monodomain model in cardiac electrophysiology. Our method is an extension of the operator splitting technique with adaptive reaction (OSTAR) where adaptive explicit integration is applied for the integration of the diffusion term too. Our approach is simple, yet efficient, and provides higher computational speed as compared to OSTAR. By introducing adaptive time stepping also for the explicit solution of the diffusion term, the conditional stability limitation can be overcome with small computational overhead, since the largest part of the required computational time is associated with the solution of the decoupled reaction term. The combination of adaptive time integration for both the diffusion and reaction terms in the DAETI method allows achieving a computational speed-up of up to x2.7 in 2D tissue simulations and x4.4 in 3D biventricular simulations when compared to simulations performed with OST without time adaptivity. Using adaptive integration only for the reaction term in the OSTAR method, the computational speed-up is limited to x1.4 and x1.8 for the same 2D and 3D simulations. Based on these findings, we can conclude that the DAETI method is notably more computationally efficient than OSTAR, especially for 3D realistic models.

The DAETI method demonstrates good numerical accuracy compared to OST without time adaptivity. In 2D tissue simulations, no differences in APD$_{90}$ calculated at the center of the tissue are found between the two methods. CV is increased by $6.4\%$ (0.0063 cm/ms in DAETI, 0.0595 cm/ms in OST). These results are as expected, since small voltage differences in the AP upstroke due to numerical approximation contribute mainly to CV and not to APD$_{90}$, which is calculated from the time point associated with maximum AP derivative and it is largely controlled by the recovery process \cite{qu1999}. The OSTAR method leads to increased CV by $3.2\%$ as compared to OST. From these findings we can conclude that the CV increase in DAETI due to the adaptive integration of the diffusion term is $3.2\%$. In \cite{qu1999}, the increase in CV due to time adaption for the reaction term is less than $2\%$ even for a time integration step as large as $dt = 0.4$ ms. In our study, the increase in CV is larger even if the time integration step $dt=0.1$ ms is smaller. It should, however, be noted that simulations in \cite{qu1999} are performed in 1D cables, while ours correspond to 2D and 3D tissues, which suggests that the increase in CV associated with adaptive integration of the reaction term may be larger for higher dimensions.

In 3D simulations using biventricular geometries, very good agreement is found between DAETI and either OSTAR or OST in terms of both LAT and APD$_{90}$. Taking OST as a reference, both DAETI and OSTAR render error values of $e_{\textrm{LAT}}$ and $e_{\textrm{APD}}$ of the order of E$^{-3}$ and E$^{-4}$, respectively, at baseline and under LBBB conditions. The activation pattern characteristic of LBBB conditions, as well as the longer time required to fully complete ventricular activation under these conditions, are well represented by DAETI, OSTAR and OST, providing results in agreement with previously published data \cite{van1978,iglesias2016}.

The findings of this study confirm that the proposed DAETI method can be effectively used to simulate cardiac electrophysiology under physiological and pathophysiological conditions, with similar numerical accuracy and remarkably higher efficiency than OSTAR, especially for 3D simulations. Importantly, by applying adaptive integration to the explicit solution of the diffusion term, the conditional stability limitation of the explicit schemes is overcome. For all the above characteristics, the DAETI method is suggested as an attractive technique for the solution of large scale problems. Since explicit schemes are highly parallelizable, we expect to observe a significantly higher efficiency gain in a parallel computing architecture implementation. In future work, we plan to implement an MPI parallel version of the DAETI method and compare its scalability against semi-implicit methods in large scale 3D applications.

\section{Conclusion} \label{sec:conclude}
A dual adaptive explicit time integration (DAETI) method is proposed to solve the monodomain model in 2D and 3D cardiac electrophysiology simulations. DAETI is based on the simple yet efficient realization that adaptive time integration can be used for the explicit solution of the diffusion term, on top of the reaction term, in the monodomain model after decoupling these two terms by the operator splitting technique. In a set of 2D and 3D simulations of cardiac electrical propagation in health and disease, the DAETI method is shown to render results of similar accuracy but improved computational efficiency than operator splitting-based methods with and without adaptive integration of the reaction term. As a fully explicit technique, the application of DAETI is expected to provide even higher efficiency gain in parallel implementations due to its straightforward parallelization and high scalability.  

\section*{Acknowledgements}
This work was supported by the European Research Council under grant agreement ERC-StG 638284, by Ministerio de Ciencia e Innovaci\'on (Spain) through project PID2019-105674RB-I00 and by European Social Fund (EU) and Arag\'on Government through BSICoS group (T39\_20R) and project LMP124-18. Computations were performed by the ICTS NANBIOSIS (HPC Unit at University of Zaragoza).

\bibliographystyle{unsrt}  
\bibliography{main}

\end{document}